\documentclass[12pt]{amsart}

\newtheorem{THM}{Theorem}[section]
\newtheorem{LMA}[THM]{Lemma}
\newtheorem{PROP}[THM]{Proposition}
\newtheorem{CORO}[THM]{Corollary}
\newtheorem{CONJ}[THM]{Conjecture}
\newtheorem{EG}[THM]{Example}

\numberwithin{equation}{section}

\usepackage{amsmath}
\usepackage{amssymb}
\usepackage{euscript}

\newcommand{\showon}{\begin{eqnarray}}
\newcommand{\showoff}{\end{eqnarray}}

\newcommand{\goesto}{\rightarrow}

\newcommand{\blam}{\boldsymbol{\lambda}}
\newcommand{\bmu}{\boldsymbol{\mu}}
\newcommand{\eps}{\varepsilon}
\newcommand{\one}{\boldsymbol{1}}

\newcommand{\A}{\EuScript{A}} 

 \newcommand{\CC}{\mathbb{C}}

\newcommand{\D}{\EuScript{D}} 
\newcommand{\E}{\EuScript{E}} \newcommand{\e}{\mathrm{e}}

\newcommand{\kB}{k_{\mathrm{B}}}

  \newcommand{\NN}{\mathbb{N}}

 \newcommand{\RR}{\mathbb{R}}
 \renewcommand{\S}{\EuScript{S}}

\renewcommand{\u}{\mathbf{u}}

 \newcommand{\x}{\mathbf{x}}

 \newcommand{\ZZ}{\mathbb{Z}}
\newcommand{\z}{\mathbf{z}}

\begin{document}

\title[Weighted enumeration of subgraphs]
{Weighted enumeration of spanning subgraphs with degree constraints}
\author{David G. Wagner}
\address{Department of Combinatorics and Optimization\\
University of Waterloo\\
Waterloo, Ontario, Canada\ \ N2L 3G1}
\email{\texttt{dgwagner@math.uwaterloo.ca}}
\thanks{Research supported by the Natural
Sciences and Engineering Research Council of Canada under
operating grant OGP0105392.}

\begin{abstract}
The Heilmann-Lieb Theorem on (univariate) matching polynomials states that
the polynomial $\sum_k m_k(G) y^k$ has only real nonpositive zeros, in which
$m_k(G)$ is the number of $k$-edge matchings of a graph $G$.  There is a stronger
multivariate version of this theorem.  We provide a general method by which
``theorems of Heilmann-Lieb type'' can be proved for a wide variety of polynomials
attached to the graph $G$.  These polynomials are multivariate generating functions
for spanning subgraphs of $G$ with certain weights and constraints imposed, and the
theorems specify regions in which these polynomials are nonvanishing.  Such theorems
have consequences for the absence of phase transitions in certain probabilistic
models for spanning subgraphs of $G$.
\end{abstract}

\maketitle

\section{Introduction.}

Let $G=(V,E)$ be a finite graph, possibly with loops or multiple edges.
For each natural number $k\in\NN$, let $m_k(G)$ denote the number of
$k$-edge matchings in $G$.  The univariate Heilmann-Lieb Theorem \cite{HL}
states that all zeros of the polynomial $\mu(G;y)=\sum_k m_k(G) y^k$ lie on the
negative real axis.  A stronger multivariate version has variables
$\x=\{x_v:\ v\in V\}$, one for each vertex, and concerns the polynomial
$$\widetilde{\mu}(G;\x)=\sum_M \x^{\deg(M)}$$
in which the sum is over all matchings $M$ of $G$,  $\deg(M):V\goesto\NN$
is the degree function of $M$, and for any $f:V\goesto\NN$
$$\x^f=\prod_{v\in V} x_v^{f(v)}.$$
The multivariate Heilmann-Lieb Theorem \cite{HL} states that if $|\arg(x_v)|<\pi/2$
for all $v\in V$ then $\widetilde{\mu}(G;\x)\neq 0$.  One sees that this
implies the univariate version by means of the relation
$$\mu(G;y)=\widetilde{\mu}(G;y^{1/2}\one)$$
(which follows from the Handshake Lemma).

The purpose of this paper is to apply some standard results from the analytic
theory of complex polynomials to provide a general method by which
``theorems of Heilmann-Lieb type'' can easily be deduced.  The
multivariate Heilmann-Lieb Theorem itself appears as the simplest
-- and prototypical -- special case of the method.  Other direct applications
provide multivariate extensions of previous results of the author \cite{Wa1},
and of results of Ruelle \cite{Ru1,Ru2}. A variety of new results also appear
as natural special cases.

In the remainder of this Introduction we
describe the general combinatorial situation we will consider.  In Section 2
we gather the necessary results from the analytic theory of complex polynomials.
In Section 3 we state and prove the main theorem of the paper.  Section 4 illustrates
this result with several applications, including the previously known examples
mentioned above.  In Section 5 we explain an interpretation of the polynomials we
consider as partition functions, by analogy with the Boltzmann-Gibbs formalism
in statistical mechanics.  Results like those in Section 4 imply that when the
thermodynamic limit of the free energy exists it must be analytic in certain regions
of the complex plane.  As noted by Lee and Yang \cite{LY,YL}, this has implications
for the absence of phase transitions in these models (which enumerate spanning
subgraphs subject to certain weights and constraints).  A more thorough
investigation of the phase structure of these models would be very interesting,
but must be left for a later paper.

It is a pleasure to thank my good friend Alan Sokal, the anonymous referee,
whose detailed positive criticism of an earlier form of this paper prompted me to
rewrite it completely.  The result is, I think, much improved.\\

The general framework we consider is that of a finite graph $G=(V,E)$ (possibly
with loops or multiple edges) and a set of weights $\blam=\{\lambda_e:\ e\in E\}$
on the edges of $G$.  These weights can for some purposes be considered as
indeterminates, but will usually be taken to be complex numbers, and often
will be nonnegative real numbers.  (In combinatorial applications it is
most natural to set all the edge-weights equal to one.)  The starting point for the
theory is the elementary identity
\showon
\Omega(G,\blam;\x)=\prod_{vew\in E}(1+\lambda_e x_v x_w)=
\sum_{H\subseteq E} \blam^H \x^{\deg(H)}.
\showoff
In this formula, the product is over the set of all edges $e\in E$, and the
notation $vew$ indicates that the ends of $e$ are the vertices $v$ and $w$
(note that $v=w$ is possible).  The sum is over the set of all spanning subgraphs
$(V,H)$ of $G$, each of which is determined by its edge-set $H\subseteq E$.
As above $\deg(H):V\goesto\NN$ is the degree function of $H$, and we use the
shorthand notations
$$\blam^H=\prod_{e\in H}\lambda_e$$
and
$$\x^{\deg(H)}=\prod_{v\in V} x_v^{\deg(H,v)}.$$

This $\Omega(G,\blam;\x)$ is a relatively structureless object, since it sums over
all spanning subgraphs without preference.  On the other hand, the product formula
allows one to make very precise statements about its zero-set (as a subset
of $\CC^V$).  To make use of this, we introduce a sequence of
\emph{activities} at each vertex $v\in V$:
\showon
\u^{(v)}=(u_0^{(v)}, u_1^{(v)},..., u_d^{(v)})\ \ \ \ (d=\deg(G,v))  
\showoff
which can be any complex numbers (usually taken to be nonnegative reals).
With these activities specified, a spanning subgraph $H\subseteq E$ will be
given the weight
\showon
\u_{\deg(H)}=\prod_{v\in V} u_{\deg(H,v)}^{(v)}
\showoff
and we will consider the correspondingly weighted version of $\Omega(G,\blam;\x)$:
\showon
Z(G,\blam,\u;\x)=\sum_{H\subseteq E}\blam^H \u_{\deg(H)}\x^{\deg(H)}.
\showoff
For example, if at every vertex we take $u_0=u_1=1$ and $u_k=0$ for all $k\geq 2$,
then 
$$
\u_{\deg(H)} = \left\{
\begin{array}{ll}
1 & \mathrm{if}\ H\ \mathrm{is\ a\ matching},\\
0 & \mathrm{otherwise},
\end{array}\right.
$$
and $Z(G,\blam,\u;\x)$ is an edge-weighted version of the multivariate matching
polynomial $\widetilde{\mu}(G;\x)$ above.

The strategy in what follows is to begin with information about the zero-set
of $\Omega(G,\blam;\x)$ and to impose conditions on the vertex activities
$\u^{(v)}$ that are sufficient to imply similar information about the zero-set
of $Z(G,\blam,\u;\x)$.  To realize this plan, we need a few results from the
analytic theory of complex polynomials.

\section{Complex polynomials.}

The technique we use is known as \emph{Schur-Szeg\H{o} composition}.
We do not make use of the most general possible result, but
for thoroughness of exposition we derive what is needed
from the Grace-Szeg\H{o}-Walsh Coincidence Theorem.  For a more
complete treatment see Sections 15 and 16 of Marden \cite{Ma} and Chapters
3 and 5 of Rahman and Schmeisser \cite{RS}.

Let $F(\z)$ be a polynomial in complex variables $\z:=\{z_{v}:\ 
v\in V\}$.  For a subset $\A\subset\CC$, we say that
$F$ is \emph{$\A$--nonvanishing} if either $F\equiv 0$, or
$z_v\in\A$ for all $v\in V$ implies that
$F(\z)\neq 0$.  In the case that $F\not\equiv 0$ we say that $F$ is
\emph{strictly} $\A$--nonvanishing.

\begin{LMA}
Let $\A$ be nonempty, connected and open.
Let $F_n(\z)$ be a sequence of strictly $\A$--nonvanishing
polynomials indexed by positive integers, and assume that
the limit $F(\z)=\lim_{n\goesto\infty} F_n(\z)$ exists.
Then $F$ is $\A$--nonvanishing.
\end{LMA}
\begin{proof}
Each $F_n$ is analytic and strictly nonvanishing on the subset $\A^V$ of $\CC^V$.
Since these functions are polynomials, the convergence to $F$ is uniform on
compact subsets of $\CC^{V}$.  By Hurwitz's Theorem (Theorem 1.3.8 of \cite{RS}),
either $F$ is identically zero or $F$ is nonvanishing on $\A^V$ as well.
\end{proof}

\begin{LMA}
Let $\A$ be nonempty, connected and open.
Let $F(\z)$ be an $\A$--nonvanishing polynomial, and let
$w\in V$.  If $z_{w}$ is fixed at a complex value $\xi$ in the closure of $\A$,
then the resulting polynomial in the variables
$\{z_{v}:\ v\in V\smallsetminus\{w\}\}$ is $\A$--nonvanishing.
\end{LMA}
\begin{proof}
The result is trivial if $F\equiv 0$, so assume instead that $F$ is strictly
$\A$--nonvanishing.  Let $(\xi_n:\ n=1,2,...)$ be a sequence with each
$\xi_n\in\A$ such that $\lim_{n\goesto \infty} \xi_n=\xi$.
Note that for all $n \geq 1$ the specialization $z_w= \xi_n$ results in
a polynomial $F_n$ that is strictly $\A$--nonvanishing
in the variables $\{z_{v}:\ v\in V\smallsetminus\{w\}\}$.
The sequence $(F_n:\ n\geq 1)$ satisfies the hypothesis of Lemma 2.1,
from which the result follows.
\end{proof}

We are concerned mostly with the following open subsets of $\CC$.\\
$\bullet$\  For $0<\theta\leq\pi$, the open sector
\showon
\S[\theta]=\{z\in\CC:\ \ z\neq 0\ \mathrm{and}\ |\arg(z)|<\theta\}
\showoff
centered on the positive real axis.  (For $z\neq 0$ we use the value of the
argument in the range $-\pi<\arg(z)\leq\pi$.)\\
$\bullet$\ For $\kappa>0$, the open interior of a disk
\showon
\kappa\D:=\{z\in\CC:\ |z|<\kappa\}.
\showoff
$\bullet$\ Also for $\kappa>0$, the open exterior of a disk
\showon
\kappa\E:=\{z\in\CC:\ |z|>\kappa\}.
\showoff
When $\kappa=1$ we more simply write just $\D$ and $\E$.

A \emph{circular region} in $\CC$ is a proper subset that is either open
or closed and is bounded by either a circle or a straight line.
A polynomial $F(\z)=F(z_1,...,z_d)$ is \emph{multiaffine} if each
variable occurs at most to the first power.  The polynomial $F(\z)$ is
\emph{symmetric} if it is invariant under every permutation of the variables.
The \emph{elementary symmetric functions} of the variables $\z=(z_1,...,z_d)$
are
\showon
e_j(\z)=\sum_{1\leq i_1< i_2<\cdots < i_j\leq d} z_{i_1} z_{i_2}\cdots z_{i_j}.
\showoff
A multiaffine symmetric polynomial $F(z_1,...,z_d)$ is thus a linear combination of
the elementary symmetric functions $e_j(\z)$ for $0\leq j\leq d$.

\begin{PROP}[Grace-Szeg\H{o}-Walsh]
Let $F(z_1,...,z_d)$ be a multiaffine symmetric polynomial, and let $\A$ be a
circular region.  Assume that either $\A$ is convex or the degree of $F$ is $d$.
Then, for any values $\zeta_1,...,\zeta_d\in\A$ there exists a value
$\zeta\in\A$ such that
$$
F(\zeta_1,...,\zeta_d)=F(\zeta,...,\zeta).
$$
\end{PROP}
For a proof in the case that $\deg F = d$, see 
Theorem 15.4 of \cite{Ma} or Theorem 3.4.1b of \cite{RS}.  The theorem also
holds when $\deg F < d$ with the additional hypothesis that $\A$ is convex,
as explained in Theorem 2.12 of \cite{COSW}.

For an elaboration of the ideas of Proposition 2.4, see Lemma 5.5.4 and Theorem
5.5.5 of \cite{RS}.
\begin{PROP}[Schur-Szeg\H{o}]
Let $P(z)=\sum_j c_j z^j$ and $K(z)=\sum_{j=0}^d \binom{d}{j} u_j z^j$
be polynomials in one complex variable $z$, with $\deg P \leq d$,
and let $Q(z)=\sum_{j=0}^d u_j c_j z^j$.\\
\textup{(a)}\ For any $0\leq\alpha<\pi/2$, if $P(z)$ is $\S[\pi/2]$-nonvanishing
and $K(z)$ is $\S[\pi-\alpha]$-nonvanishing, then $Q(z)$ is
$\S[\pi/2-\alpha]$-nonvanishing.\\
\textup{(b)}\ For any $\kappa>0$ and $\rho>0$, if $P(x)$ is $\rho\D$-nonvanishing
and $K(z)$ is $\kappa\D$-nonvanishing, then $Q(z)$ is $\kappa\rho\D$-nonvanishing.\\
\textup{(c)}\ For any $\kappa>0$ and $\rho>0$, if $P(x)$ is $\rho\E$-nonvanishing
and $K(z)$ is $\kappa\E$-nonvanishing and $\deg K = d$,
then $Q(z)$ is $\kappa\rho\E$-nonvanishing.
\end{PROP}
\begin{proof}
The conclusions are trivial if $Q\equiv 0$, so we may assume that $Q\not\equiv 0$.

We begin by proving part (a) in the case that $K(0)\neq 0$.
In this case we have
\showon
K(z)=C\prod_{i=1}^d(1+\theta_i z)
\showoff
for some complex numbers $C\neq 0$ and $\theta_1,...,\theta_d$ such that
either $\theta_i=0$ or $|\arg(\theta_i)|\leq
\alpha$ for each $1\leq i\leq d$.  Consider the \emph{$d$-th polarization
of $P(z)$}:\ this is the multiaffine symmetric polynomial
$\widetilde{P}(\z)=\widetilde{P}(z_1,...,z_d)$ obtained
from $P(z)$ by replacing each monomial $z^j$ by the normalized $j$-th elementary
symmetric function $\binom{d}{j}^{-1}e_j(\z)$.  Since $\deg P \leq d$,
it follows that
\showon
\widetilde{P}(z,z,...,z)=P(z)
\showoff
as polynomials in $z$.  Since $P(z)$ is $\S[\pi/2]$-nonvanishing and
$\S[\pi/2]$ is a circular region, it follows from (2.6) and Proposition 2.4
that $\widetilde{P}(\z)$ is also $\S[\pi/2]$-nonvanishing.
Now, consider complex numbers $\zeta_1,...,\zeta_d\in\S[\pi/2-\alpha]$.
For each $1\leq i\leq d$, either $\theta_i\zeta_i=0$ for all
$\zeta_i\in\S[\pi/2-\alpha]$ or $|\arg(\theta_i\zeta_i)|<\pi/2$ for all
$\zeta_i\in\S[\pi/2-\alpha]$.  From Lemma 2.2,
it follows that if $\widetilde{P}(\theta_1 z_1,...,\theta_d z_d)\not\equiv 0$
then $\widetilde{P}(\theta_1\zeta_1,...,\theta_d\zeta_d)\neq 0$
for every choice of $\zeta_1,...,\zeta_d\in\S[\pi/2-\alpha]$.  That is, it
follows that $\widetilde{P}(\theta_1 z_1,...,\theta_d z_d)$ is
$\S[\pi/2-\alpha]$-nonvanishing.  A short calculation using the fact
that $\binom{d}{j}u_j= C e_j(\theta_1,...,\theta_d)$ verifies that
\showon
Q(z) = C \widetilde{P}(\theta_1 z,...,\theta_d z),
\showoff
and therefore $Q(z)$ is $\S[\pi/2-\alpha]$-nonvanishing, as desired.

To handle the case in which $K(0)=0$, let $r$ be the multiplicity of $0$
as a root of $K(z)$ and write
\showon
K(z)=C z^r \prod_{i=1}^{d-r}(1+\theta_i z).
\showoff
For a positive integer $N$ let
\showon
K_N(z)=C N^{-r}(1+Nz)^r \prod_{i=1}^{d-r}(1+\theta_i z).
\showoff
and let $Q_N(z)$ be the polynomial in the conclusion constructed from
$P(z)$ and $K_N(z)$.  By the case we have done already, each $Q_N(z)$ is
$\S[\pi/2-\alpha]$-nonvanishing.  Taking the limit as $N\goesto\infty$,
Lemma 2.1 implies that $Q(z)$ itself is also $\S[\pi/2-\alpha]$-nonvanishing.

The proof of part (b) is similar.  Since $K(z)$ is $\kappa\D$-nonvanishing
we have $K(0)\neq 0$, and so we can write $K(z)$ as in equation (2.5) with all
$|\theta_i|\leq 1/\kappa$.  Again we consider the $d$-th polarization 
$\widetilde{P}(\z)$ of $P(z)$.  Since $P(z)$ is $\rho\D$-nonvanishing and
$\rho\D$ is a circular region, Proposition 2.3 and equation (2.6) imply
that $\widetilde{P}(\z)$ is $\rho\D$-nonvanishing.  It follows that
$\widetilde{P}(\theta_1 z_1,...,\theta_d z_d)$ is $\kappa\rho\D$-nonvanishing,
and from equation (2.7) we conclude that $Q(z)$ is $\kappa\rho\D$-nonvanishing,
as desired.

The proof of part (c) repeats the same pattern once more.
Begin with $K(z)$ expressed as in equation (2.8) -- since  $K(z)$ is
$\kappa\E$-nonvanishing, each $|\theta_i|\geq 1/\kappa$.  We work with
the polynomials $K_N(z)$ defined in equation (2.9) with $N\geq 1/\kappa$.
Since $P(z)$ is $\rho\E$-nonvanishing and
$\rho\E$ is a circular region and $\deg \widetilde{P}=d$, Proposition 2.3 and
equation (2.6) imply that $\widetilde{P}(\z)$ is $\rho\E$-nonvanishing. 
It follows that
$$\widetilde{P}(\theta_1 z_1,...,\theta_{d-r} z_{d-r}, Nz_{d-r+1},...,Nz_d)$$
is $\kappa\rho\E$-nonvanishing,
and from equation (2.7) we conclude that $Q_N(z)$ is $\kappa\rho\E$-nonvanishing.
Taking the limit as $N\goesto\infty$ (using Lemma 2.1) we conclude that $Q(z)$ is
$\kappa\rho\E$-nonvanishing, as desired.
\end{proof}
The polynomial $Q(z)$ in the conclusion of Proposition 2.4 is the
\emph{Schur-Szeg\H{o} composition} of $P(z)$ and $K(z)$.

\section{The main result.}

Consider a graph $G=(V,E)$ with complex edge weights $\blam$.
We begin with some easy information about the zero-set of the polynomial
$\Omega(G,\blam;\x)$ defined in equation (1.1).
\begin{PROP}
Let $G=(V,E)$ be a graph with complex edge weights $\blam$.\\
\textup{(a)}\  If $\lambda_e\geq 0$ for each $e\in E$ then
$\Omega(G,\blam;\x)$ is $\S[\pi/2]$-nonvanishing.\\
\textup{(b)}\  If $|\lambda_e|\leq\lambda_{\max}$ for each $e\in E$ then
$\Omega(G,\blam;\x)$ is $\lambda_{\max}^{-1/2}\D$-nonvanishing.\\
\textup{(c)}\  If $|\lambda_e|\geq\lambda_{\min}$ for each $e\in E$ then
$\Omega(G,\blam;\x)$ is $\lambda_{\min}^{-1/2}\E$-nonvanishing.
\end{PROP}
\begin{proof}
In each case, each factor $1+\lambda_e z_v z_w$ in the product form for
$\Omega(G,\blam;\z)$ is seen to be nonvanishing in the appropriate region,
from which the result follows.
\end{proof}

Now assume that we also have a sequence of activities $\u^{(v)}$
at each vertex $v\in V$, as in equation (1.2).  The information about
these activities that we will use is recorded in the set of
\emph{key polynomials}
\showon
K_v(z)=\sum_{j=0}^d \binom{d}{j} u_j^{(v)} z^j
\showoff
in which $d=\deg(G,v)$.  There is one key polynomial for each vertex $v\in V$.

\begin{THM}
Let $G=(V,E)$ be a graph, with complex edge weights $\blam$, and with vertex
activities $\u$ encoded by the key polynomials $K_v(z)$ ($v\in V$).\\
\textup{(a)}  Fix $0\leq \alpha<\pi/2$.
If $\lambda_e\geq 0$ for each $e\in E$ and $K_v(z)$ is
$\S[\pi-\alpha]$-nonvanishing for each $v\in V$,
then $Z(G,\blam,\u;\x)$ is $\S[\pi/2-\alpha]$-nonvanishing.\\
\textup{(b)}  Fix $\kappa>0$ and $\lambda_{\max}>0$.
If  $|\lambda_e|\leq\lambda_{\max}$ for each $e\in E$ and $K_v(z)$ is
$\kappa\D$-nonvanishing for each $v\in V$,
then $Z(G,\blam,\u;\x)$ is $(\kappa/\lambda_{\max}^{1/2})\D$-nonvanishing.\\
\textup{(c)}  Fix $\kappa>0$ and $\lambda_{\min}>0$.
If  $|\lambda_e|\geq\lambda_{\min}$ for each $e\in E$ and $K_v(z)$ is
$\kappa\E$-nonvanishing and $\deg K_v(z)=\deg(G,v)$ for each $v\in V$,
then $Z(G,\blam,\u;\x)$ is $(\kappa/\lambda_{\min}^{1/2})\E$-nonvanishing.
\end{THM}
\begin{proof}
Identify the vertices $V$ with the numbers $V=\{1,2,...,n\}$ arbitrarily.
Define a sequence of polynomials $F_0(\x)$, $F_1(\x)$,..., $F_n(\x)$ as follows.
$F_0(\x)=\Omega(G,\blam;\x)$, and for all $1\leq v\leq n$, $F_v(\x)$ is the
Schur-Szeg\H{o} composition of $F_{v-1}(\x)$ regarded as a polynomial
in the variable $x_v$ (the other variables being absorbed into the coefficients)
with $K_v(x_v)$.  One sees by induction that for $0\leq r\leq n$:
\showon
F_r(\x)=\sum_{H\subseteq E}\blam^H
\left(\prod_{v=1}^r u_{\deg(H,v)}^{(v)}\right)\x^{\deg(H)},
\showoff
so that $F_n(\x)=Z(G,\blam,\u;\x)$.

We give the details to finish the proof of part (a) -- the arguments for parts
(b) and (c) are completely analogous.  We prove by induction on $1\leq v\leq n$
that if $(\zeta_j:\ 1\leq j\leq n)$ are complex numbers such that\\
$\bullet$\ $|\arg(\zeta_j)|<\pi/2-\alpha$ for all $1\leq j<v$, and\\
$\bullet$\ $|\arg(\zeta_j)|<\pi/2$ for all $v<j\leq n$,\\
then
\showon
F_{v-1}(\zeta_1,...,\zeta_{v-1},x_v,\zeta_{v+1},...,\zeta_n)
\showoff
is $\S[\pi/2]$-nonvanishing.  The basis of induction follows from Proposition 3.1(a)
and Lemma 2.2.  The induction step follows from Proposition 2.4(a) and Lemma 2.2.
Finally, from the statement that whenever all $\zeta_i\in\S[\pi/2-\alpha]$,
then $F_{n-1}(\zeta_1,...,\zeta_{n-1},x_n)$ is
$\S[\pi/2]$-nonvanishing, we conclude by one more application of Proposition 2.4(a)
that $F_{n}(\x)$ is $\S[\pi/2-\alpha]$-nonvanishing, as desired.
\end{proof}

The univariate specialization of Theorem 3.2 is an important consquence.
\begin{CORO}
Adopt the notation of Theorem $3.2$.\\
\textup{(a)}\ Under the hypotheses of Theorem \textup{3.2(a)},\\
$Z(G,\blam,\u; y^{1/2}\one)$ is $\S[\pi-2\alpha]$-nonvanishing.\\
\textup{(b)}\ Under the hypotheses of Theorem \textup{3.2(b)},\\
$Z(G,\blam,\u; y^{1/2}\one)$ is $(\kappa^2/\lambda_{\max})\D$-nonvanishing.\\
\textup{(c)}\ Under the hypotheses of Theorem \textup{3.2(c)},\\
$Z(G,\blam,\u; y^{1/2}\one)$ is $(\kappa^2/\lambda_{\min})\E$-nonvanishing.
\end{CORO}

\section{Applications.}

Throughout this section, consider a graph $G=(V,E)$ with complex edge
weights $\blam$ and vertex activities $\u$ encoded by the key polynomials
$K_v(z)$ ($v\in V$).

\begin{EG}[Heilmann-Lieb \cite{HL}]\emph{
Assume that all edge weights are nonnegative reals, and that at each vertex
$u_0=u_1=1$ and $u_k=0$ for all $k\geq 2$.  The key polynomial at a 
vertex of degree $d$ in $G$ is $K_v(z)= 1+dz$, which is $\S[\pi]$-nonvanishing.
Theorem 3.2(a) (with $\alpha=0$) implies that $Z(G,\blam,\u;\x)$ is $\S[\pi/2]$-nonvanishing -- this is the multivariate Heilmann-Lieb theorem.
Corollary 3.3(a) implies that $Z(G,\blam,\u;y^{1/2}\one)$ is
$\S[\pi]$-nonvanishing -- this is the univariate Heilmann-Lieb theorem.
}\end{EG}

\begin{EG}[Wagner \cite{Wa1}]\emph
{
Assume that all edge weights are nonnegative reals, and that two functions
$f,g:V\goesto\NN$ are given such that $f(v)\leq g(v)\leq f(v)+1$ for each
$v\in V$.  Fix the vertex activities to be
\showon
u_k^{(v)}=\left\{
\begin{array}{ll}
1 & \mathrm{if}\ f(v)\leq k\leq g(v),\\
0 & \mathrm{otherwise}.
\end{array}\right.
\showoff
As in Example 4.1, each key $K_v(z)$ is $\S[\pi]$-nonvanishing.
Theorem 3.2(a) (with $\alpha=0$) implies that $Z(G,\blam,\u;\x)$ is
$\S[\pi/2]$-nonvanishing -- this result is new.
Corollary 3.3(a) implies that $Z(G,\blam,\u;y^{1/2}\one)$ is
$\S[\pi]$-nonvanishing -- when $\blam\equiv\one$ this is
Theorem 3.3 of \cite{Wa1}.
}\end{EG}

\begin{EG}[Ruelle \cite{Ru1,Ru2}]\emph{
Assume that all edge weights are nonnegative reals, and that at each vertex
$u_0=u_2=1$, $u_1=u$, and $u_k=0$ for all $k\geq 2$.  The key polynomial at a 
vertex of degree $d$ in $G$ is $K_v(z)= 1+d u z+ \binom{d}{2}z^2$.
For $d\geq 2$, the zeros of this polynomial are at
$$z_{\pm}=\frac{-2}{d-1}\left(u\pm\sqrt{u^2-2+2/d}\right).$$
\indent
When $u=1$, all the keys $K_v(z)$ are $\S[3\pi/4]$-nonvanishing, and
Theorem 3.2(a) (with $\alpha=\pi/4$) implies that $Z(G,\blam,\u;\x)$ is
$\S[\pi/4]$-nonvanishing -- this is new.  
Corollary 3.3(a) implies that $Z(G,\blam,\u;y^{1/2}\one)$ is
$\S[\pi/2]$-nonvanishing -- when $\blam\equiv\one$
this is a slight weakening of Proposition 1 of \cite{Ru1}.\\ \indent
If $G$ has maximum degree $\Delta$ and $u\geq\sqrt{2-2/\Delta}$,
then all the keys $K_v(z)$ are $\S[\pi]$-nonvanishing, and
Theorem 3.2(a) (with $\alpha=0$) implies that $Z(G,\blam,\u;\x)$ is
$\S[\pi/2]$-nonvanishing -- this result is new.  
Corollary 3.3(a) implies that $Z(G,\blam,\u;y^{1/2}\one)$ is
$\S[\pi]$-nonvanishing -- when $\blam\equiv\one$
this is Proposition 2 of \cite{Ru1}.\\ \indent
Ruelle's method produces more detailed information than ours, but only for
particular choices of the vertex activities.  A systematic extension of his
method that handles all the cases we consider would be very interesting.
}\end{EG}

\begin{EG}\emph{
Assume that all edge weights are nonnegative reals, and that two functions
$f,g:V\goesto\NN$ are given such that $f(v)\leq g(v)\leq f(v)+2$ for each
$v\in V$.  Fix the vertex activities as in equation (4.1).
Then each key $K_v(z)$ is $\S[2\pi/3]$-nonvanishing.
Theorem 3.2(a) (with $\alpha=\pi/3$) implies that $Z(G,\blam,\u;\x)$ is
$\S[\pi/6]$-nonvanishing, and Corollary 3.3(a) implies that
$Z(G,\blam,\u;y^{1/2}\one)$ is $\S[\pi/3]$-nonvanishing. 
}\end{EG}

\begin{EG}\emph{
Assume that all edge weights are nonnegative reals, and that two functions
$f,g:V\goesto\NN$ are given such that $f(v)\leq g(v)\leq f(v)+3$ for each
$v\in V$.  Fix the vertex activities as in equation (4.1).
If every vertex of $G$ has degree at most $\Delta$ then there is a small
angle $\eps>0$ such that each key $K_v(z)$ is $\S[\pi/2+\eps]$-nonvanishing.
To see this, the keys with at most three terms pose no problems (by Examples
4.2 and 4.4).  A key with four terms has the form
$$
K(z)=\binom{d}{f} z^f + \binom{d}{f+1} z^{f+1} + \binom{d}{f+2} z^{f+2}
+ \binom{d}{f+3} z^{f+3},
$$
and the inequality
$$\binom{d}{f+1}\binom{d}{f+2} > \binom{d}{f}\binom{d}{f+3}$$
ensures that the only zero of $K(z)$ with nonnegative real part is at the
origin.  Since $\Delta$ is fixed, only finitely many key polynomials need to
be considered -- taking the smallest positive argument of the (nonzero) zeros
of these to be $\pi/2+\eps$ gives the desired angle.\\ \indent
Theorem 3.2(a) implies that $Z(G,\blam,\u;\x)$ is
$\S[\eps]$-nonvanishing, and Corollary 3.3(a) implies that
$Z(G,\blam,\u;y^{1/2}\one)$ is $\S[2\eps]$-nonvanishing.
}\end{EG}

\begin{EG}\emph{
In Examples 4.1 and 4.2 we concluded that the polynomial
$Z(y)=Z(G,\blam,\u;y^{1/2}\one)$ had only real (and nonpositive) zeros.
Let $N_j=N_j(G,\blam,\u)$ be the coefficient of $y^j$ in this polynomial.
It is ``folklore'' that if $Z(y)$ is $\S[2\pi/3]$-nonvanishing, then
$$N_i N_k\neq 0 \ \ \mathrm{implies\ that}\ \ N_j\neq 0\ \ \mathrm{for\ all}\
\ i\leq j\leq k$$
and
$$N_j^2\geq N_{j+1}N_{j-1}\ \ \mathrm{for\ all}\ \ j.$$
This property (\emph{logarithmic concavity with no internal zeros})
is very useful for obtaining good approximations to the sequence
$(N_j)$ (see \cite{FK,K,HL}, for example).\\ \indent
If all the keys $K_v(z)$ are $\S[5\pi/6]$-nonvanishing then $Z(y)$ is
$\S[2\pi/3]$-nonvanishing.  However, this hypothesis on the keys is
unreasonably strong.  Consider a key of the form
$$K(z)=\binom{d}{j-1}z^{j-1}+\binom{d}{j}z^{j}+\binom{d}{j+1}z^{j+1}$$
with $1\leq j\leq d-1$,
corresponding to three consecutive permissible degrees.  A short calculation
shows that this is $\S[5\pi/6]$-nonvanishing if and only if $2j(d-j)\leq d+2$.
This happens only for the pairs $(j,d)$ with $d\leq 4$ and $j=1$ or $j=d-1$.\\
\indent
Nonetheless, I venture the following conjecture.
}\end{EG}

\begin{CONJ}
Let $G=(V,E)$ be a finite graph, and let $f,g:V\goesto\NN$ be any two functions.
Fix the vertex activities $\u$ as in $(4.1)$.  Then the sequence of coefficients
$(N_j)$ of $Z(G,\one,\u;y^{1/2}\one)$ is logarithmically concave with no
internal zeros.
\end{CONJ}

\begin{EG}\emph{
Assume that all the edge weights have unit modulus, that $G$ is
$2k$-regular, and that the key at each vertex is
\showon
K(z) = 1 + \binom{2k}{k} z^k + u z^{2k}.
\showoff
If $4 u \geq \binom{2k}{k}^2$ then every zero of $K(z)$ has modulus
$\kappa=u^{-1/2k}$.  Parts (b) and (c) of Theorem 3.2 imply that
$Z(G,\blam,\u;\x)$ is both $\kappa\D$- and $\kappa\E$-nonvanishing.
Corollary 3.3 implies that every zero of $Z(G,\blam,\u;y^{1/2}\one)$
has modulus $u^{-1/k}$.
}\end{EG}

\begin{EG}\emph{
Assume that all the edge weights have unit modulus, and that
$\deg K_v(z)=\deg(G,v)$ and every zero of $K_v(z)$ has unit modulus,
for each vertex $v\in V$.  Parts (b) and (c)
of Theorem 3.2 imply that $Z(G,\blam,\u;\x)$ is both $\D$-
and $\E$-nonvanishing.  Corollary 3.3 implies that every zero of
$Z(G,\blam,\u;y^{1/2}\one)$ has unit modulus.\\ \indent
In particular, these hypotheses evidently hold if $\blam\equiv\one$ and
the key polynomials are given by $K_v(z)=1+z+z^2+\cdots+z^{\deg(G,v)}$
for each $v\in V$.  Thus we conclude that every zero of
$$\sum_{H\subseteq E}\frac{y^{\#H}}{\prod_{v\in V}\binom{\deg(G,v)}{\deg(H,v)}}$$
has unit modulus.
}\end{EG}

\section{Analogy with statistical mechanics.}

We conclude with an interpretation of $Z(G,\blam,\u; y^{1/2}\one)$
inspired by analogy with the (canonical ensemble) partition functions
in statistical mechanics.
For simplicity, we restrict attention to a graph $G=(V,E)$ that is $d$-regular,
in which the edge weights $\blam\equiv\one$ are all one and the
activities are the same at every vertex (that is, all the key polynomials are equal).
The extension to the general case is straightforward.
 
The ``configuration space'' is the set of all spanning subgraphs
of $G$.  The \emph{energy} $U(H)$ of a spanning subgraph $H\subseteq E$ depends
on $d+2$ real parameters $J$ and $\bmu=(\mu_0, \mu_1,..., \mu_d)$, as follows:
\showon
U(H)=J\cdot\#H+\sum_{j=0}^d \mu_j\cdot\#V_j(H),
\showoff
in which $V_j(H)$ is the set of vertices of degree $j$ in $H$.
The quasi-physical interpretation of this is that $J$ is the energy of a single
edge, and $\mu_j$ is the ``chemical potential'' energy of a vertex of degree $j$.
With $T>0$ denoting absolute temperature, and $\beta=1/\kB T$ where
$\kB$ is Boltzmann's constant, the \emph{Boltzmann weight} of $H$ is
$$\e^{-\beta U(H)}$$ and the \emph{partition function} is
\showon
Z_G(\beta,J,\bmu)=\sum_{H\subseteq E} \e^{-\beta U(H)}.
\showoff

This can be interpreted as defining a family  of probability measures
(parameterized by $\beta$, $J$, and $\bmu$) on
the set of all spanning subgraphs of $G$:\  a spanning subgraph $H\subseteq E$
is chosen at random with probability $\e^{-\beta U(H)}/Z_G(\beta,J,\bmu)$.
A short computation shows that, for $H$ chosen according to this distribution,
the expected number of edges is
\showon
\langle \#H \rangle =
-\frac{1}{\beta}\frac{\partial}{\partial J}\log Z_G(\beta,J,\bmu)
\showoff
and the expected number of vertices of degree $j$ is
\showon
\langle \#V_j(H) \rangle =
-\frac{1}{\beta}\frac{\partial}{\partial \mu_j}\log Z_G(\beta,J,\bmu)
\showoff

To continue with the analogy we consider a sequence of graphs
$G_1$, $G_2$,... that converges to an infinite, locally finite,
limit graph $\Gamma$. (The precise definition of convergence is not
important for this discussion -- the prototypical example is that,
as $n\goesto\infty$, the Cartesian product $C_n^r$ of $r$ cycles of length
$n$ should converge to the infinite graph $\ZZ^r$ with edges of Euclidean
length one.)  We will further assume that the
``thermodynamic limit'' \emph{(Helmholtz) free energy}
\showon
f_\Gamma(\beta,J,\bmu)=-\frac{1}{\beta}
\lim_{n\goesto\infty} \frac{1}{\#V(G_n)}\log  Z_{G_n}(\beta,J,\bmu)
\showoff
exists.  As in the Lee-Yang theory \cite{LY,YL},
points in the parameter space at which the
free energy fails to be analytic can be interpreted as phase transitions
between differing qualitative properties of a random spanning subgraph of 
$\Gamma$.  From the form of (5.5) we see that $f_\Gamma$ can fail to be
analytic only at an accumulation point of the union of the zero-sets of
all the $Z_{G_n}(\beta,J,\bmu)$ ($n\geq 1$).  From the probabilistic
interpretation, we are most interested in such accumulation points for which
all the parameters $(\beta,J,\bmu)$ are real.

The partition functions can be expressed as polynomials in the variables
\showon
y=\e^{-\beta J}\ \ \ \mathrm{and}\ \ \ u_j=\e^{-\beta\mu_j}\ (0\leq j\leq d).
\showoff
In fact, a tiny calculation shows that in these variables
\showon
Z_G(\beta,J,\bmu)=Z(G,\one,\u;y^{1/2}\one)
\showoff
with the RHS as defined in (1.4).  The point $y=1$ corresponds to
$\beta J=0$, which is the infinite-temperature limit.  If $J>0$ then
$y=0$ is the zero-temperature limit, and if $J<0$ then $y\goesto+\infty$
is the zero-temperature limit.  The positive real axis is thus the ``physically''
relevant part of the complex $y$-plane.  If all the chemical potentials $\mu_j$
are real then all the activities $u_j$ are positive reals.  A zero activity
$u_j=0$ corresponds to an infinite chemical potential $\mu_j=+\infty$,
which means that a vertex of degree $j$ is forbidden.  Notice that the
activity $u_j=\e^{-\beta\mu_j}$ also depends on temperature except when
$\mu_j$ is $+\infty$ or $0$:\ this is the case precisely when $u_j\in\{0,1\}$.

In this context, Corollary 3.3 has the following immediate consequence,
the proof of which is omitted.
\begin{PROP}
Let $(G_n:\ n\geq 1)$ be a sequence of $d$-regular graphs, and let
$\beta>0$ and $J\in\RR$ and $\bmu\in\RR^{d+1}$ be such that the limit
$(5.5)$ exists.  Form the key polynomial
$$K(z)=K(\beta,\bmu;z)=\sum_{j=0}^d \binom{d}{j} u_j z^j$$
with $(u_j)$ as in $(5.6)$.\\
\textup{(a)}\ If there exists $\eps>0$ such that $K(z)$
is $\S[\pi/2+\eps]$-nonvanishing then
$f_\Gamma$ is analytic at $(\beta,J,\bmu)$ for all $J\in\RR$.\\
\textup{(b)}\ If $\kappa>0$ is such that $K(z)$ is $\kappa\D$-nonvanishing
then $f_\Gamma$ is analytic at $(\beta,J,\bmu)$ for all
$$J>-\frac{2}{\beta}\log\kappa.$$
\textup{(c)}\ If $\kappa>0$ is such that $K(z)$ is $\kappa\E$-nonvanishing
and of degree $d$ then $f_\Gamma$ is analytic at $(\beta,J,\bmu)$ for all
$$J<-\frac{2}{\beta}\log\kappa.$$
\end{PROP}

Finally, we revisit some of the examples of Section 4, maintaining as well
the assumptions of Proposition 5.1.

\begin{EG}\emph{
With the key polynomial $K(z)$ as in Example 4.3,
let $u=\e^{-\beta\mu}$.  If $\mu<+\infty$ (that is, if $u>0$) then $K(z)$ is
$\S[\pi/2+\eps]$-nonvanishing for all $\beta\geq 0$, so that $f_\Gamma$ is
analytic at $(\beta,J,\bmu)$ for all $J\in\RR$.  In this case there is
no phase transition at any nonzero temperature.
On the other hand, if $\mu=+\infty$ (that is, if $u=0$) then both zeros
of $K(z)$ have modulus $\kappa=\binom{d}{2}^{-1/2}$, so that
$f_\Gamma$ is analytic at $(\beta,J,\bmu)$ for all
$$J>\frac{1}{\beta}\log\binom{d}{2}.$$
In this case there is no phase transition provided that the
temperature $T$ is sufficiently low compared to the edge energy $J$.
}\end{EG}

\begin{EG}\emph{
With the key polynomial $K(z)$ as in Example 4.5,
$K(z)$ is $\S[\pi/2+\eps]$-nonvanishing for all $\beta\geq 0$, so
that $f_\Gamma$ is analytic at $(\beta,J,\bmu)$ for all $J\in\RR$. 
Thus, there is never a phase transition in this model.
}\end{EG}

\begin{EG}\emph{
With the key polynomial $K(z)$ as in Example 4.8, let $u=\e^{-\beta\mu}$
and $d=2k$.
If $\binom{2k}{k}^2\leq 4u$ then all the zeros of $K(z)$ have modulus $\kappa=u^{-1/2k}$
and $K(z)$ has degree $d$.  The only point on the positive $y$-axis at which
$f_\Gamma$ could fail to be analytic is at $y=u^{-1/k}$.  In terms of the ``physical''
parameters, this says that if
\showon
-\beta\mu > 2\log\binom{2k}{k}-\log 4
\showoff
then a phase transition can occur only at $J=-\mu/k$.  The inequality (5.8) requires
that $\mu<0$ (so that vertices of degree $2k$ in $H$ are energetically favoured) and
that $\beta$ is sufficiently large (so that the temperature is sufficiently low).
If this is the case then a phase transition can occur only when the edge energy $J$
and chemical potential $\mu$ are tuned to satisfy $J=-\mu/k$.
}\end{EG}

\begin{EG}\emph{
With the key polynomial $K(z)$ as in Example 4.9,
all the zeros of $K(z)$ have modulus one and $K(z)$ has degree $d$.
The only point on the positive $y$-axis at which
$f_\Gamma$ could fail to be analytic is at $y=1$.  In terms of the ``physical''
parameters, this says that a phase transition can occur only at $\beta J=0$ --
that is, only in the infinite temperature limit.
}\end{EG}

As these examples illustrate, Proposition 5.1 sees very little about the
limit graph $\Gamma$ -- in fact, only the degree of $\Gamma$ is relevant.
(On the other hand, the existence of the limit $f_\Gamma$ does depend on
the structure of $\Gamma$.)
Thus, for example, Proposition 5.1 can not tell the difference between the
3d cubical lattice and the 2d triangular grid -- both graphs are regular
of degree six.  Of course, in truth one expects that for any given model,
the free energies of these two graphs will have different phase diagrams.
Accounting for more detailed structural properties of $\Gamma$ remains an
interesting open problem.

\end{document}